\documentclass[11pt]{article}
\usepackage{amsmath,amsthm}

\def\refeq#1{\if\workingver y(\ref{#1})-[[#1]]\else(\ref{#1})\fi}
\def\refth#1{\if\workingver y\ref{#1}-[[#1]]\else\ref{#1}\fi}
\def\mylabel#1{\if\workingver y\label{#1}{\bf\ \ [[#1]]\ \ }
\else\label{#1}\fi}
\def\mybibitem#1{\if\workingver y\bibitem{#1}{\bf\ \ [[#1]]\ \ }
\else\bibitem{#1}\fi}

\def\institute#1{\gdef\@institute{#1}}

\def\institutename{\par
 \begingroup
 \parskip=\z@
 \parindent=\z@
 \setcounter{@inst}{1}%
 \def\and{\par\stepcounter{@inst}%
 \noindent$^{\the@inst}$\enspace\ignorespaces}%
 \setbox0=\vbox{\def\thanks##1{}\@institute}%
 \ifnum\c@@inst=1\relax
 \else
   \setcounter{footnote}{\c@@inst}%
   \setcounter{@inst}{1}%
   \noindent$^{\the@inst}$\enspace
 \fi
 \ignorespaces
 \@institute\par
 \endgroup}

\newtheorem{thm}{Theorem}
\newtheorem{lem}[thm]{Lemma}

\newtheorem{cor}[thm]{Corollary}
\newtheorem{prop}[thm]{Proposition}

\newtheorem{rem}[thm]{Remark}
\newtheorem{conj}[thm]{Conjecture}

\def\begeq#1{\begin{equation}\mylabel{#1}}
\def\endeq{\end{equation}}

\def\begalg{\begin{alg}}
\def\endalg{\end{alg}}

\let\workingver=n

\newcommand{\beq}{\begin{equation}}

\def\ga{\gamma}
\def\de{\delta}
\def\al{\alpha}
\def\be{\beta}

\begin{document}

\title{Matrix Powers of Column-Justified Pascal Triangles and Fibonacci
Sequences}
\author{Rhodes Peele, Pantelimon St\u anic\u a\\
\small Auburn University Montgomery, Department of Mathematics\\
\small Montgomery, AL 36117, USA\\
\small e-mail: \em peelhor@mindspring.com, stanpan@strudel.aum.edu}
\date{\today}
\maketitle

\baselineskip=1.9\baselineskip

\section{Motivation}

It is known that, if $L_n$, respectively $R_n$ are
the $n\times n$ matrices
with the $(i,j)$th entry the binomial coefficient
$\binom{i-1}{j-1}$, respectively $\binom{i-1}{n-j}$, then
$L_n^2\equiv I_n\pmod{2}$, respectively $R_n^3\equiv I_n
\pmod{2}$, where $I_n$ is the identity matrix of dimension $n>1$
(see for instance the problem P10735 in
the May 1999 issue of American Mathematical Monthly).

The entries of $L_n$ form a left-justified Pascal's Triangle and
the entries
of $R_n$ result from taking the mirror-image of this Triangle
with respect to its first column.

The questions that we ask are: can this result be extended to
other primes or better yet,
is it possible to find a closed form for the entries of powers
of $L_n$ and $R_n$?

$L_n$ succumbs easily as we shall see in our first result. $R_n$ in
turn fights back,
since closed forms for its powers are not found. However, we show a
beautiful connection between
matrices similar to $R_n$ and the Fibonacci numbers. If $n=2$, the
 connection is easily seen, since
\[
\begin{split}
R_2^e &=\left(
\begin{matrix}
F_{e-1} & F_{e}\\
F_{e} & F_{e+1}
\end{matrix}
\right).
\end{split}
\]
A simple consequence of our results is that
 the order of $L_n$ modulo a prime $p$ is $p$,
and the order of $R_n$ modulo $p$ divides $4$ times the
entry point of the Fibonacci sequence modulo~ $p$.

\section{Higher Powers of $L_n$ and $R_n$}

The first approach that comes to mind is to find a closed form
for all entries of
powers of $L_n,R_n$.
It is not  difficult to obtain all the powers of $L_n$.
Denoting the entries of the $e$-th power of $L_n$ by $l_{i,j}^{(e)}$
we can prove
\begin{thm}
\label{higher_L}
The entries of $L_n^e$ are
\begin{equation}
\label{l_higher}
l_{i,j}^{(e)}=e^{i-j}\binom{i-1}{j-1}.
\end{equation}
\end{thm}
\begin{proof}
We use induction on $e$. The result is certainly true for $e=1$.
Now, using induction and matrix multiplication,
\begin{eqnarray*}
l_{i,j}^{(e+1)}
&=& \sum_{s=1}^n \binom{i-1}{s-1} e^{s-j} \binom{s-1}{j-1}\\
&\stackrel{\cite{R},\ p.\, 3}{=} &
\sum_{s=1}^n e^{s-j} \binom{i-1}{j-1}\binom{i-j}{i-s}\\
&=& \binom{i-1}{j-1}\sum_{k=0}^{i-j} e^{i-j-k} \binom{i-j}{k}=
\binom{i-1}{j-1}(e+1)^{i-j}.
\end{eqnarray*}
\end{proof}

To prove a similar result for $R_n$ is no easy matter. In fact,
except for a few
lower-dimensional cases and a few of its rows/columns, simple
closed forms for the
entries of $R_n^e$ are not found.


 In the sequel we consider the tableau with entries
$a_{ij},\, i\geq 1,\, j\geq 0$ satisfying
\begin{equation}
\label{cond_matrix}
a_{i,j-1}=a_{i-1,j-1}+a_{i-1,j},
\end{equation}
with boundary conditions $a_{1,n}= 1,\, a_{1,j}=0,\, j\neq n$. We
shall use the following consequences of the boundary conditions
and recurrence \refeq{cond_matrix}: $a_{i,j}=0$ for $i+j\leq n$,
and $a_{i,n+1}=0, 1\leq i\leq n$ (in fact we use {\em only} these
consequences and \refeq{cond_matrix}). The matrix $R$ will be
defined as $\left(a_{i,j}\right)_{i=1...n,j=1...n}$. We treat the
second and third powers first, since it gives us the idea about
the general case.
To clear up the mysteries of some of the steps in
our calculations we will refer
to matrix multiplication as {\em m.m.} and
boundary conditions as {\em b.c.}
\begin{lem}
\label{lem2}
The entries of the matrix $R^2$ satisfy
\begin{equation}
\label{second_power}
b_{i,j+1}=b_{i-1,j+1}+2 b_{i-1,j}-b_{i,j},\ 2\leq i\leq n,\
1\leq j\leq n-1
\end{equation}
and the entries of $R^3$ satisfy
\begin{equation}
\label{third_power}
c_{i+1,j}=2c_{i,j}+3 c_{i,j-1}-2 c_{i+1,j-1},\ 1\leq i\leq n-1,\
2\leq j\leq n.
\end{equation}
\end{lem}
\begin{proof}
Using matrix multiplication and \refeq{cond_matrix} we obtain
\begin{equation}
\begin{split}
\label{eq7}
b_{i,j+1}
\stackrel{m.m.}{=}& \sum_{s=1}^n a_{i,s} a_{s,j+1}
\stackrel{\refeq{cond_matrix}}{=}\sum_{s=1}^n a_{i,s}
(a_{s+1,j}-a_{s,j})=\\
=& \sum_{s=1}^n a_{i,s} a_{s+1,j} -\sum_{s=1}^n a_{i,s} a_{s,j}
\stackrel{m.m.}{=}
\sum_{s=1}^n a_{i,s} a_{s+1,j}-b_{i,j}.
\end{split}
\end{equation}
Therefore, denoting $S_{i,j}=\sum_{s=1}^n a_{i,s}a_{s+1,j}$, we obtain
\beq
\label{eq6}
b_{i,j+1}+b_{i,j}=S_{i,j}.
\end{equation}
If $2\leq i\leq n$ and $1\leq j\leq n$,
\begin{eqnarray*}
S_{i,j}
&=&
\sum_{s=1}^n \left( a_{i-1,s}+a_{i-1,s+1} \right) a_{s+1,j}
\stackrel{t=s+1}{=}
S_{i-1,j}+
\sum_{t=2}^{n+1} a_{i-1,t}a_{t,j}\\
&\stackrel{m.m.}{=}&
S_{i-1,j}+b_{i-1,j}+a_{i-1,n+1}a_{n+1,j}-a_{i-1,1}a_{1,j}
\stackrel{b.c.}{=}
S_{i-1,j}+ b_{i-1,j}.
\end{eqnarray*}
Using \refeq{eq6} in the previous recurrence, we obtain
\[
b_{i,j+1}+b_{i,j}=b_{i-1,j+1}+b_{i-1,j}+b_{i-1,j},
\]
which gives us \refeq{second_power}.

If the relations \refeq{second_power} are satisfied,
we obtain, for $j\geq 2$,
\[
\begin{split}
c_{i,j}
\stackrel{m.m.}{=}& \sum^n_{s=1} a_{i,s} b_{s,j}
\stackrel{\refeq{second_power}}{=}\sum_{s=1}^n a_{i,s}(b_{s-1,j}+2
b_{s-1,j-1}-b_{s,j-1})\\
=& \sum_{s=1}^n a_{i,s} b_{s-1,j}+2\sum_{s=1}^n  a_{i,s}
b_{s-1,j-1}-c_{i,j-1}\\
=&\, T_{i,j}+2T_{i,j-1}-c_{i,j-1},
\end{split}
\]
where $T_{i,j}=\sum^n_{s=1} a_{i,s} b_{s-1,j}$.
Furthermore, for $i\leq n-1$,
\begin{eqnarray*}
T_{i,j}
&\stackrel{\refeq{cond_matrix}}{=}&
\sum_{s=1}^n \left( a_{i+1,s-1}-a_{i,s-1} \right) b_{s-1,j}\\
&\stackrel{m.m.}{=}&
c_{i+1,j} +a_{i+1,0}b_{0,j}-a_{i+1,n}b_{n,j}-c_{i,j}-a_{i,0}b_{0,j}
+a_{i,n}b_{n,j}\\
&\stackrel{\refeq{cond_matrix}}{=}&
c_{i+1,j}-c_{i,j}+a_{i,1}b_{0,j}-a_{i,n+1}b_{n,j}
\stackrel{b.c.}{=} c_{i+1,j}-c_{i,j}.
\end{eqnarray*}
Therefore,
\begin{eqnarray*}
c_{i,j}&=&T_{i,j}+2T_{i,j-1}-c_{i,j-1}\\
&=&c_{i+1,j}-c_{i,j}+2c_{i+1,j-1}-2c_{i,j-1}-c_{i,j-1},
\end{eqnarray*}
which will produce the equations \refeq{third_power}.
\end{proof}

\begin{cor}
\label{Cor2_3}
The entries of the second and third power of
$R$ can be expressed in terms of the
entries of the previous row:
\[
\begin{split}
b_{i+1,j} & = b_{i,j}-\sum_{k=1}^{j-1}  (-1)^k b_{i,j-k},\\
c_{i+1,j} & = 2c_{i,j}+\sum_{k=1}^{j-1} (-1)^k 2^{k-1} c_{i,j-k}.
\end{split}
\]
\end{cor}

We have wondered if relations similar to \refeq{second_power}
or \refeq{third_power}
are true for higher powers of $R$. It turns out that
\begin{thm}
\label{high_pow}
The entries $a_{i,j}^{(e)}$ of the $e$-th power of $R$
satisfy the relation
\[
F_{e-1} a_{i,j}^{(e)}=F_{e} a_{i-1,j}^{(e)}+
F_{e+1} a_{i-1,j-1}^{(e)}-F_{e} a_{i,j-1}^{(e)},
\]
where $F_e$ is the Fibonacci sequence.
\end{thm}
\begin{proof}
We show first that the entries of $R^e$ satisfy a relation of the form
\begin{equation}
\label{eq77}
\delta_e a_{i,j}^{(e)}=\alpha_e a_{i-1,j}^{(e)}+\beta_e a_{i-1,j-1}^{(e)}+
\gamma_e a_{i,j-1}^{(e)}
\end{equation}
and then will proceed to find these coefficients. From Lemma
\refth{lem2} we observe that
$\delta_1=0,\alpha_1=1,\beta_1=1,\gamma_1=-1$,
$\delta_2=1,\alpha_2=1,\beta_2=2,\gamma_2=-1$ and
$\delta_3=1,\alpha_3=2,\beta_3=3,\gamma_3=-2$. Now, the
coefficients of $R^e$ satisfy, for $i,j\geq 2$,
\begin{equation}
\begin{split}
\de_{e-1}a_{i-1,j}^{(e)}
\stackrel{m.m.}{=}&
\sum_{s=1}^n \de_{e-1} a_{i-1,s}
a_{s,j}^{(e-1)}\\
\stackrel{\refeq{eq77}}{=}&
\sum_{s=1}^n a_{i-1,s}  \left(\al_{e-1}
a_{s-1,j}^{(e-1)}+\be_{e-1}a_{s-1,j-1}^{(e-1)}+\ga_{e-1}
a_{s,j-1}^{(e-1)}\right)\\
\stackrel{m.m.}{=} &
\al_{e-1} U_{i-1,j}+\be_{e-1} U_{i-1,j-1}+\ga_{e-1} a_{i-1,j-1}^{(e)},
\end{split}
\end{equation}
where $U_{i,j}=\sum_{s=1}^n a_{i,s}  a_{s-1,j}^{(e-1)}$. We evaluate, for
$2\leq i\leq n$,
\begin{eqnarray*}
 U_{i-1,j}
&\stackrel{\refeq{cond_matrix}}{=}&
\sum_{s=1}^n \left(a_{i,s-1} -a_{i-1,s-1}  \right)a_{s-1,j}^{(e-1)}\\
&\stackrel{m.m.}{=}&
a_{i,j}^{(e)}-a_{i-1,j}^{(e)} +a_{i,0}  a_{0,j}^{(e-1)}-
a_{i-1,0}a_{0,j}^{(e-1)}-a_{i,n} a_{n,j}^{(e-1)}+
a_{i-1,n} a_{n,j}^{(e-1)}\\
&\stackrel{\refeq{cond_matrix}}{=}&
a_{i,j}^{(e)}-a_{i-1,j}^{(e)}+a_{i-1,1}  a_{0,j}^{(e-1)}-a_{i-1,n+1}
a_{n,j}^{(e-1)}=a_{i,j}^{(e)}-a_{i-1,j}^{(e)},
\end{eqnarray*}
since $a_{i-1,1}=0,i\leq n$ and $a_{i-1,n+1}=0.$
Thus,
\[
\al_{e-1} a_{i,j}^{(e)}=
(\de_{e-1}+\al_{e-1})a_{i-1,j}^{(e)}+(\be_{e-1}-\ga_{e-1})
a_{i-1,j-1}^{(e)}-
\be_{e-1}a_{i,j-1}^{(e)}.
\]
 Therefore, we obtain the following system of sequences
\begin{equation}
\begin{split}
\de_e =\, & \al_{e-1}\\
\al_e=\, & \al_{e-1}+\de_{e-1}\\
\be_e=\,& \be_{e-1}-\ga_{e-1}\\
\ga_e=\, & - \be_{e-1}.
\end{split}
\end{equation}
From this we deduce $\de_{e}=F_{e-1},\al_{e}=F_{e},\be_e=F_{e+1},
\ga_e=-F_{e}$,
where $F_e$ is the Fibonacci sequence, with $F_0=0,F_1=1$.
\end{proof}

Corollary \refth{Cor2_3} can be generalized,
with a little more work and anticipating \refeq{first_row},
to obtain the elements in
the $(i+1)$th row of $R^e$,
in terms of the elements in the previous row.
\begin{prop}
We have
\[
F_{e-1} a_{i+1,j}^{(e)}=F_ea_{i,j}^{(e)}-
\sum_{k=1}^{j-1} (-1)^{k+e} \frac{F_e^{k-1}}{F_{e-1}^k}a_{i,j-k}^{(e)}.
\]
\end{prop}

\section{Higher Powers of $L_n$ and $R_n$ modulo a prime $p$}

As before, let $L_n$, respectively, $R_n$ be defined as the matrices with
entries $\binom{i-1}{j-1}$, respectively, $\binom{i-1}{n-j}$. We
use the notation ``$matrix \equiv a\pmod p$" with the meaning
``$matrix \equiv a I_n\pmod p$".

We ask the question of whether the order of $L_n, R_n$ modulo a
prime $p$ is finite.
We can prove easily a result for $L_n$ using Theorem \refth{higher_L}.
\begin{thm}
The order of $L_n$ $(n\geq 2)$ modulo  $p$ is $p$.
\end{thm}
\begin{proof}
We showed that the entries of $L_n^e$ are $l_{i,j}^{(e)}=e^{i-j}
\binom{i-1}{j-1}$,
for any integer $e$. Thus, the entries on the principal diagonal
of $L_n^e$ are all 1. If $i\neq j$, then $p\,| l_{i,j}^{(p)}$.
Assume that there is an integer $e$ with $0<e<p$, such that
$p\,| l_{i,j}^{(e)}$, for all $i\not= j$.
Take $i=2,j=1$. Thus, $p\,| e$, a contradiction.
Therefore, the integer $p$  is
the least integer $e> 0$ for which $p\,| l_{i,j}^{(e)}$,
for all $i\not= j$,
 which proves our assertion.
\end{proof}

We can prove the finiteness of the order of $R_n$ modulo $p$  in a simple
manner. By the Pigeonhole Principle, there exist $s<t$, such that
 $R_n^s\equiv R_n^t\pmod {p}$.
Since $R_n$ is an invertible matrix
($\det {R_n}=(-1)^n\not\equiv 0\pmod p$),
 $R_n^{t-s}\equiv I_n\pmod p$.  More precise results will be proved next.
In order to do that we need some known facts about the period of
the Fibonacci sequence.
It was shown that the period of the Fibonacci sequence modulo
$m$ (not necessarily prime) is less than or equal to $6m$ (with equality
holding for infinitely many values of $m$) (P. Freyd, American
Mathematical Monthly, December 1990, E 3410, with a solution
provided in the March 1992 issue of AMM).
In the case of a prime, the result can be strengthened (see the
next theorem).
The least integer $n\not =0$
with the property $m\,|F_n$ is called the {\em entry point} modulo $m$.

\def\P{{\cal P}}

In \cite{Bloom} and \cite{Wall} the authors obtain (see also
\cite{Vajda1989}, Ch. VI-VII, for a more updated source)
\begin{thm}[Bloom-Wall]
Denote the period of the Fibonacci sequence modulo $p$ by $\P(p)$.
Let $p$ be an odd prime with $p\not= 5$. If $p\equiv \pm
1\pmod {5}$, then the period $\P(p)\, |\, (p-1)$. If $p\equiv \pm
3\pmod {5}$, then the entry point $e\,| (p+1)$ and the period
$\P(p)\, |\, 2(p+1)$.
\end{thm}

\begin{rem}
For $p=2$ the entry point and the period is $3$.
In the case $p=5$ the entry point is $5$ and the period is $20$.
\end{rem}

\begin{thm}
\label{Th10}
\begin{sloppypar}
If $e$ is the entry point modulo $p$ of $F_e$, then
$R_{2k}^e\equiv (-1)^{(k+1)e} F_{e-1}I_{2k}\pmod p$ and
$R_{2k+1}^e\equiv (-1)^{ke}I_{2k+1} \pmod p$.
Moreover, $R_n^{4e}\equiv I_n\pmod{p}$.
\end{sloppypar}
\end{thm}

\begin{proof}
We prove by induction on $e$ that the elements on the first row
and first column
of $R_n^e$ are
\begin{equation}
\label{first_row}
\begin{split}
a_{1,j}^{(e)} &=\binom{n-1}{j-1} F_{e-1}^{n-j} F_e^{j-1}\\
a_{i,1}^{(e)} &=F_{e-1}^{n-i} F_e^{i-1}.
\end{split}
\end{equation}
First we deal with the elements in the first row. The top
 equation is certainly true for $e=1$,
if we define $0^0=1$.
 Now,
\begin{eqnarray*}
a_{1,j}^{(e+1)}
&\stackrel{m.m.}{=}&
\sum_{s=1}^n a_{1,s}^{(e)} a_{s,j}=\sum_{s=1}^n F_{e-1}^{n-s}
F_e^{s-1}\binom{n-1}{s-1}\binom{s-1}{n-j}\\
&=&
\sum_{s=1}^n  F_{e-1}^{n-s} F_e^{s-1} \binom{n-1}{j-1}\binom{j-1}{n-s}\\
&=&
F_e^{n-1} \binom{n-1}{j-1}
\sum_{s=1}^n \left( \frac{F_{e-1}}{F_e}\right)^{n-s} \binom{j-1}{n-s}\\
&=&
F_e^{n-1} \binom{n-1}{j-1} \left(1+\frac{F_{e-1}}{F_e}\right)^{j-1}=
F_{e}^{n-j} F_{e+1}^{j-1} \binom{n-1}{j-1}.
\end{eqnarray*}
Again by induction we prove the result for the elements in the
first column. The case $e=1$ is checked easily.  Then,
\begin{eqnarray*}
a_{i,1}^{(e+1)}
&\stackrel{m.m.}{=}&
\sum_{s=1}^n a_{i,s} a_{s,1}^{(e)}=\sum_{s=1}^n
\binom{i-1}{n-s} F_{e-1}^{n-s} F_e^{s-1} \\
&=&
F_e^{n-1} \sum_{s=1}^n \left(\frac{F_{e-1}}{F_e}
\right)^{n-s}\binom{i-1}{n-s}=F_e^{n-1}
\left(1+\frac{F_{e-1}}{F_e}\right)^{i-1}\\
&=&
F_e^{n-i} F_{e+1}^{i-1}.
\end{eqnarray*}

Let $e$ be the entry point modulo $p$ of the Fibonacci sequence.
By Bloom-Wall's result, $e\leq p+1$.
Using Theorem \refth{high_pow}, we obtain
\[
F_{e-1} a_{i,j}^{(e)}\equiv F_{e} a_{i-1,j}^{(e)}+
F_{e+1} a_{i-1,j-1}^{(e)}-F_{e} a_{i,j-1}^{(e)}\pmod{p}.
\]
Thus,
\begin{equation}
\label{propag}
F_{e-1} a_{i,j}^{(e)}\equiv F_{e+1} a_{i-1,j-1}^{(e)}\pmod{p}.
\end{equation}
Since $F_{e-1}+ F_{e}=F_{e+1}$, $p\,|F_e$ and $p\not| F_{e-1}$, we obtain
$F_{e-1}\equiv F_{e+1}\pmod p$ and
\begin{equation}
\label{propag2}
a_{i,j}^{(e)}\equiv a_{i-1,j-1}^{(e)}\pmod p.
\end{equation}
We see, from what was proved above that modulo $p$,
the elements on the first row and column of $R_n\pmod p$ are all zero,
except for the one on the first position, which is
$F_{e-1}^{n-1}\not\equiv 0\pmod p$.
Using
\refeq{propag2} we get $R_n^e\equiv F_{e-1}^{n-1} I_n \pmod p$.
Using Cassini's identity
$F_{e-1}F_{e+1}-F_e^2=(-1)^e$ (see \cite{GKP}, p. 292),
we obtain
$F_{e-1}^2\equiv F_{e+1}^2\equiv (-1)^e\pmod p$.
If $n=2k$, then
\[
F_{e-1}^{n-1}
=F_{e-1}^{2k-1}\equiv \left( F_{e-1}^2\right)^k F_{e-1}^{-1}
\equiv  (-1)^{ke} F_{e-1}^{-1}\equiv (-1)^{(k+1)e} F_{e-1}\pmod p.
\]
If $n=2k+1$, then \[
F_{e-1}^{n-1}
= F_{e-1}^{2k}\equiv \left( F_{e-1}^2\right)^k
\equiv  (-1)^{ke} \pmod p.
\]
The previous two congruences replaced in
$R_n^{e}\equiv F_{e-1}^{n-1} I_n\pmod p$, will
give the first two assertions of our theorem.

It is well known (a very particular case of Matijasevich's lemma) that
\[
F_{2e-1}=F_{e-1}^2+F_e^2\equiv F_{e-1}^2\pmod {F_{e}^2},
\]
so $F_{2e-1}^2\equiv 1\pmod p$.
Thus, since $F_m$ divides $F_{sm}$ for all $m,s$
(in particular for $s=2,\, m=e$),
it follows that
$F_{2e}\equiv 0\pmod p$ and
\[
R_n^{4e}=\left(R_n^{2e}\right)^2\equiv \left(F_{2e-1}^2\right)^{n-1}
I_n\equiv I_n \pmod p.
\]
\end{proof}

\begin{rem}
We remark here the fact that the bound $4e$ for the order of $R$ is tight.
That can be seen by taking, for example, the prime $13$, since
the entry point for the
Fibonacci sequence is $7$, and the order of $R_{4k}$ is $28$.
\end{rem}

Using some Elementary Number Theory we can prove
\begin{thm}
 If  $p\,|F_{p-1}$,
 then $R_n^{p-1} \equiv I_n \pmod p$.
\end{thm}
\begin{proof}
We observe that, since $p\, |F_{p-1}$, we have
 $p\equiv \pm 1\pmod {5}$, otherwise $p\equiv \pm 2\pmod {5}$, and
 by Bloom-Wall's Theorem
 the entry point $e$ divides $p+1$. Therefore, $e\,| p-1$ and
 $e\,| p+1$. Therefore, $e$ must be 2. That is not possible,
 since $F_2=1$, which is not divisible by any prime.
So,   $p\equiv \pm 1\pmod {5}$
and
$F_p\equiv F_{p-2}\pmod p$. Thus
\[
R_n^{p-1}\equiv F_{p-2}^{n-1} I_n \equiv F_p^{n-1} I_n\pmod p.
\]
By the previous Bloom-Wall's Theorem, $\P(p)\,|(p-1)$, therefore
$F_{p-1}\equiv 0$, $F_{p}\equiv 1$,
$F_{p+1}\equiv 1$, etc. Thus
\[
R_n^{p-1}\equiv F_p^{n-1} I_n\equiv I_n\pmod p.
\]
\end{proof}

Another interesting result is the following
\begin{thm}
If  $p\,|F_{p+1}$,
then $R_{2k+1}^{p+1} \equiv I_{2k+1} \pmod p$ and
$R_{2k}^{p+1} \equiv -I_{2k} \pmod p$.
\end{thm}
\begin{proof}
Let $p=2$. The entry point of the Fibonacci sequence modulo $2$ is $e=3$.
Since $F_2=1$,
Theorem \refth{Th10} shows the result in this case.
 Assume $p>2$.
We know that in this case we must have $p\equiv \pm 2\pmod {5}$.
Using the known formula (see for instance Hardy and Wright \cite{HW},
Theorem 180)
\[
F_j=2^{1-j}\left[ \binom{j}{1}+
5\binom{j}{3}+5^2\binom{j}{5}+\cdots\right],
\]
 taking $j=p$, and using Fermat's Little  Theorem,
 $2^{p-1}\equiv 1\pmod p$, we obtain
\[
F_p\equiv 5^{(p-1)/2}\binom{p}{p}\equiv -1\pmod p,
\]
since for the primes $\equiv \pm 2\pmod 5$, $5$ is a quadratic
nonresidue.

When $n$ is odd,
\[
R_n^{p+1}\equiv F_p^{n-1} I_n\equiv
(F_p^2)^{\frac{n-1}{2}}I_n\equiv I_n\pmod p.
\]
Consider the case of $n$  even. Since
$F_p\equiv -1\pmod p$, we have
\[
R_n^{p+1}\equiv F_p^{n-1} I_n \equiv (-1)^{n-1} I_n\equiv -I_n\pmod p.
\]
\end{proof}
The proofs of the previous two theorems imply
\begin{cor}
If $p\equiv \pm 1\pmod{5}$ and $p-1$ is the entry point for the Fibonacci
sequence modulo $p$,
then the period  is exactly $p-1$.
If $p\equiv \pm 2\pmod{5}$ and $p+1$ is the entry point for the Fibonacci
sequence modulo $p$,
then the period  is exactly $2(p+1)$.
\end{cor}

\begin{cor}
The order of $R_n\pmod p$ is less than or equal to $2(p+1)$ and
the bound is met.
\end{cor}
\begin{proof}
If $p\equiv \pm 1\pmod {5}$, then the order of $R_n\pmod p$ is
$\leq p-1$.
If $p\equiv \pm 2\pmod {5}$, then
$F_{p}\equiv -1\pmod p$. Therefore,
$R_n^{p+1}\equiv F_p^{n-1} I_n\equiv (-1)^{n-1} I_n\pmod p$. Thus
$R_n^{2p+2}\equiv I_n\pmod p$. The bound is met for all primes
$p\equiv \pm 2\pmod {5}$ and all even integers $n$.
\end{proof}

\section{Further Problems and Results}

The inverses of $R_n$ and $L_n$ are not difficult to find. We have
\begin{thm}
The inverse of
\[
L_n=\left(\binom{i-1}{j-1}\right)_{1\leq i,j\leq n}\
\text{ is}\  L_n^{-1}=\left((-1)^{i+j}\binom{i-1}{j-1}
\right)_{1\leq i,j\leq n}.
\]
The inverse of
\[
R_n=\left(\binom{i-1}{n-j}\right)_{1\leq i,j\leq n}\  \text{is}\ R_n^{-1}=
\left((-1)^{n+i+j+1}\binom{n-i}{j-1}\right)_{1\leq i,j\leq n}.
\]
\end{thm}
\begin{proof}
We have
\[
\begin{split}
& \sum_{s}^n (-1)^{i+s}\binom{i-1}{s-1}\binom{s-1}{j-1}\\
&\stackrel{\cite{R},\ p.\, 3}{=}
 \sum_{s}^n (-1)^{i+s}  \binom{i-1}{j-1}\binom{i-j}{i-s}\\
&\stackrel{k=i-s}{=}
\binom{i-1}{j-1}
\sum_{k=0}^{i-j} (-1)^{k}  \binom{i-j}{k},
\end{split}
\]
which is $0$, unless $i=j$, in which case it is $1$.
A similar analysis for $R_n$ will produce  its inverse.
\end{proof}

\def\bphi{\bar\phi}

Another approach to find a closed form of all
entries of powers of $R_n$  would be to find all eigenvalues of $R_n$, and
use the diagonalization of the matrix to find the entries of
$R_n$. We found empirically and we state as a conjecture
\begin{conj}
Denote $\phi=\frac{1+\sqrt{5}}{2},\bar\phi=\frac{1-\sqrt{5}}{2}$.
The eigenvalues of $R_n$ are:\\
1. $\displaystyle\{ (-1)^{k+i} \phi^{2i-1},
(-1)^{k+i} \bar\phi^{2i-1}  \}_{i=1,\ldots,k}$, if $n=2k$.\\
2. $\displaystyle\{ (-1)^k \}\cup\{ (-1)^{k+i} \phi^{2i},(-1)^{k+i}
\bar\phi^{2i}  \}_{i=1,\ldots,k}$, if $n=2k+1$.
\end{conj}

Another venue of research would be to study the
matrices associated to other interesting sequences: Lucas, Pell, etc., and
we will approach this matter elsewhere.

\vspace{.5cm}

\noindent{\bf Acknowledgements.}
{
The authors would like to thank the anonymous referee for his
helpful comments, which improved significantly the presentation
of the paper.
}

\noindent AMS Classification Numbers: 05A10, 11B39, 11B65, 11C20, 15A33

\end{document}